\documentclass[10pt,twoside]{article}
\usepackage{latexsym}
\def\thebibliography#1{\center{\bf\normalsize References}\list
 {[\arabic{enumi}]}{\settowidth\labelwidth{[#1]}\leftmargin\labelwidth
 \advance\leftmargin\labelsep
 \usecounter{enumi}}
 \def\newblock{\hskip .11em plus .33em minus .07em}
 \sloppy\clubpenalty4000\widowpenalty4000
 \sfcode`\.=1000\relax}
\def\tt{\hspace{-0.19cm}{\bf .}\hspace{0.19cm}}

\usepackage{graphicx}
\newcommand{\beq}[2]{\begin{equation}\label{#1}#2\end{equation}}

\newcommand{\set}[1]{\left\{#1\right\}}
\newcommand{\ip}[1]{\left\langle #1\right\rangle}
\newcommand{\ia}[1]{\left| #1\right|}
\newcommand{\thb}[1]{{\rm (#1)}}
\newcommand{\doo}[2]{\frac{d {#1}}{d {#2}}}

\def\eqif{\, {\rm if}\,\, }

\def\eqin{ \, {\rm in } \,\, }
\def\eqae{ \, {\rm a.e. } \, \,}
\def\eqst{ \, {\rm s.t. } \, \,}

\def\all{  \, \forall \, }
\def\tr{\,{\rm tr}\,}

\def\qq{\qquad}
\def\q{\quad}
\def\nnb{\nonumber}
\def\ds{\displaystyle}
\def\cd{\cdot}

\def\ol{\overline}

\def\a{\alpha}
\def\b{\beta}
\def\g{\gamma}
\def\d{\delta}

\def\ve{\varepsilon}
\def\s{\sigma}
\def\t{\tau}
\def\o{\omega}
\def\z{\zeta}
\def\G{\Gamma}
\def\O{\Omega}
\def\Th{\Theta}

\def\Bf{{\bf f}}

\def\dbR{{\mathop{\rm l\negthinspace R}}}

\def\bu{{\bar u}}
\def\bx{{\bar x}}

\def\bU{{\ol{U}}}
\def\bW{{\ol{W}}}
\def\bpsi{{\ol{\psi}}}
\def\bGO{{\ol{\O}}}
\def\bGP{{\ol{\Phi}}}

\def\mt{{\hat t}}
\def\hu{{\hat u}}
\def\hx{{\hat x}}

\def\tiu{{\tilde u}}

\def\cT{{\mathcal T}}
\def\cU{{\mathcal U}}
\def\cV{{\mathcal V}}
\def\cM{{\mathcal M}}

\newcommand{\eqref}[1]{$(\ref{#1})$}
\newcommand{\refeq}[1]{$(\ref{#1})$}

\def\Proof{\noindent{\bf Proof.\hspace{1em}}}

\def\endpf{\hfill$\Box$\vspace{0.4cm}\\ }
\def\eqif{\, {\rm if}\,\, }

\def\eqin{ \, {\rm in } \,\, }
\def\eqae{ \, {\rm a.e. } \, \,}
\def\eqst{ \, {\rm s.t. } \, \,}

\begin{document}
\title{
{\bf Second-Order Necessary/Sufficient  Conditions for Optimal
Control Problems  in the Absence of Linear Structure
\thanks{The author is supported by  NSFC (No. 10831007),
 FANEDD (No. 200522) and 973 Program (No. 2011CB808002).
}}
 }
\author{ Hongwei Lou \thanks{School of Mathematical Sciences, and LMNS, Fudan
University, Shanghai 200433, China (Email:
\texttt{hwlou@fudan.edu.cn})}
 }

\date{}
\maketitle
\begin{quote}
\footnotesize {\bf Abstract.} Second-order necessary
conditions for optimal control problems are considered, where the
``second-order" is in the sense of that Pontryagin's maximum
principle is viewed as a first-order necessary optimality condition. A
sufficient  condition for a local minimizer is also given.


\end{quote}
\normalsize

\newtheorem{Definition}{Definition}[section]
\newtheorem{Theorem}[Definition]{Theorem}
\newtheorem{Lemma}[Definition]{Lemma}
\newtheorem{Corollary}[Definition]{Corollary}
\newtheorem{Proposition}[Definition]{Proposition}
\newtheorem{Remark}{Remark}[section]
\newtheorem{Example}{Example}
\newfont{\Bbb}{msbm10 scaled\magstephalf}
\newfont{\frak}{eufm10 scaled\magstephalf}
\newfont{\sfr}{msbm7 scaled\magstephalf}

\def\theequation{1.\arabic{equation}}
\setcounter{equation}{0} \setcounter{section}{1}
\setcounter{Definition}{0} \setcounter{Remark}{0}\textbf{1.
Introduction.}  We will give a new kind of second-order necessary/sufficient
conditions for optimal controls. Second-order
necessary/sufficient conditions for optimal control
problems have been studied for a long time. There are many relevant
works. Among them, we mention the following works and the references
therein: \cite{Bon-Her1}--- \cite{Cas-Tr-Un},
\cite{G-K}---\cite{Kre}, \cite{Mal}---\cite{Wach}  and
\cite{Wang-He}. There are different definitions of ``first-order necessary/sufficient conditions"
and ``second-order necessary/sufficient conditions".
To our best knowledge, the corresponding
first-order necessary conditions to these second-order conditions in the literature
are not  Pontryagin's maximum principle. In addition, the control domains
considered there are domains (or closed domains) in $\dbR^m$. In other words, second-order necessary/sufficient
conditions for optimal controls in the literature are mainly used
to distinguish optimal controls from other singular controls in the
classical sense, not from other singular controls in
the sense of Pontryagin's maximum principle.
For the definitions of singular controls in the
classical sense and singular controls in
the sense of Pontryagin's maximum principle, see
Definitions $1$ and $2$ in \cite{G-K}, see also
\eqref{ES002}  and  \eqref{E402A}.

Before we focus on our problems, we recall the results about
necessary conditions for  minimizers of functions.

Let us consider a minimizer $x_0$ of a smooth function $f(\cd)$ on
$\bGO$, where $\bGO$ is the closure of a domain $\O\subseteq \dbR^n$.
We call a unit vector $\ell$ is an admissible direction if there
exists a $\d>0$ such that $x_0+s\ell \in \bGO$ for any $s\in
[0,\d]$. If $\ell$ is admissible, we have the following  first-order
necessary condition:
\beq{E101}{ %
  0  \leq   \lim_{s\to 0^+} {f(x_0+s\ell)-f(x_0)\over
s}=\ip{\nabla f(x_0),\ell}.
 }
When
\beq{E102}{ %
 \ip{\nabla f(x_0),\ell}=0
 }
holds, i.e., \eqref{E101} degenerates, then we can get further the
second-order necessary condition:
\beq{E103}{ %
  0  \leq   \lim_{s\to 0^+} {f(x_0+s\ell)-f(x_0)\over
s^2}={1\over 2}\ip{D^2 f(x_0)\ell,\ell},
 }
where $D^2f$ is the Hessian matrix of $f$. \if{
$$
D^2f\equiv \pmatrix{\ddp 2 f {x_1} & {\pa^2 f\over \pa x_1\pa x_2}  & \ldots &  {\pa^2 f\over \pa x_1\pa x_n}\cr
                    {\pa^2 f\over \pa x_2\pa x_1} & \ddp 2 f {x_2} & \ldots & {\pa^2 f\over \pa x_2\pa x_n}\cr
                     \vdots & \vdots & \ddots & \vdots\cr
                     {\pa^2 f\over \pa x_n\pa x_1} & {\pa^2 f\over \pa x_n\pa x_2}& \ldots &  \ddp 2 f {x_n}}.
$$}\fi
If \eqref{E102} does not hold, that is
$$
\ip{\nabla f(x_0),\ell}>0,
$$
then  \eqref{E103} does not necessarily hold .

From the above observations, we see that to yield second-order
conditions of a minimizer,  linear structure of independent
variables is needed and second-order conditions only appear when
first-order conditions degenerate.

For an optimal control problem, usually the control domain $U$ need
not have linear structure. Thus, the space $\cU_{ad}$ of control
functions need not have linear structure. Pontryagin's maximum
principle is a kind of necessary conditions that a  minimizer
satisfies.  Many people look it as the first-order necessary
condition. However,
 Pontryagin's maximum principle could not be obtained
\textbf{directly} in a way like \eqref{E102}. First, for an optimal
control $\bu(\cd)$, there is probably no ``admissible  direction"
$v(\cd)$ such that $\bu(\cd)+s v(\cd)$ is still in $\cU_{ad}$.
Secondly, even if ``admissible direction $v(\cd)$" exists, what we
could get from
$$
0\leq \lim_{s\to 0^+}{J(\bu(\cd)+sv(\cd))-J(\bu(\cd))\over s}
$$
is only a corollary of Pontryagin's maximum principle which looks like  \eqref{ES001},
where we denote $J(\cd)$ the cost functional of the optimal control
problem.

When linear structure lacks, could we replace the ``admissible
direction" by ``admissible path"?  In other words, could we replace
$\bu(\cd)+s v(\cd)$ by $u_s(\cd)\in \cU_{ad}$, which is continuous
in some sense in $s\in [0,1]$? Certainly, we can do that. Yet,
``admissible path" will immediately puzzles us on what are
first-order conditions and second-order conditions. To see this, let
us consider the function $f(\cd)$ and its minimizer $x_0$ again. Let
$\ell$ be an admissible direction such that \eqref{E102} holds. Then
choosing $x(s)=x_0+ \sqrt s \ell$, we have
\beq{E104}{ %
  0  \leq   \lim_{s\to 0^+} {f(x(s))-f(x_0)\over
s}={1\over 2}\ip{D^2 f(x_0)\ell,\ell}.
 }
Then, should we call \eqref{E104} a first-order condition?
Therefore, we think it is not a good idea to replace ``admissible
direction" by ``admissible path". In this paper, we will transform
the original optimal control problem to a new problem, which is in
fact the locally relaxed problem of the original problem. In this
new problem, the corresponding space of control functions has linear
structure and we can yield Pontryagin's maximum principle like
\eqref{E101} under this linear structure.  Then we can further
yield second-order conditions based on Pontryagin's maximum
principle.

To reveal our idea clearly, we consider simply optimal control
problems governed  by ordinary differential equations.

The rest of the paper is organized as follows: In Section 2, we will
give a  method to linearize the control space near the optimal control.
In Section 3, We will give a new proof of Pontryagin's maximum principle.
 Section 4 will be devoted to second-order necessary conditions of optimality.
 Finally, a sufficient condition for a control being a local minimizer will be given in Section 5.

\bigskip

\def\theequation{2.\arabic{equation}}
\setcounter{equation}{0} \setcounter{section}{2}
\setcounter{Definition}{0} \setcounter{Remark}{0}\textbf{2. Local
Linearization of Optimal Control Problems.} In this section, we will
linearize locally an optimal control problem along its minimizer.
Let us consider the following controlled system:
\begin{equation}\label{E201}
\left\{\begin{array}{ll}\ds \dot x(t)=f(t,x(t),u(t)),& \eqin [0,T],
\\
x(0)=x_0 &
\end{array}\right.
\end{equation}
and the following cost functional
\begin{equation}\label{E202}
J(u(\cd))=\int^T_0 f^0(t,x(t),u(t))\, dt,
\end{equation}
where $T>0$, and $u(\cd)\in \cU_{ad}$ with
\begin{equation}\label{E203}
\cU_{ad}=\set { v:[0,T]\to U\big|\, v(\cd) \q \mbox{measurable}\,}.
\end{equation}
We pose the following assumptions:

 (S1) The metric
space  $(U,\rho)$ is  separable.

 (S2) Functions  $\Bf=\pmatrix{f^0\cr f}=\pmatrix{f^0 & f^1 & f^2 & \ldots & f^n}^\top:[0,T]\times \dbR^n\times U\to\dbR^{n+1}$ are
  measurable in $t$, continuous in  $(x,u)$ and  continuously differentiable in $x$, where $B^\top$ denotes the transposition of a matrix $B$.  Moreover,  there exists a constant $L>0$ such that
\begin{equation}\label{E204}
\left\{\begin{array}{l} \ds |\Bf(t,x,u)-\Bf (t,\hx, u)|\leq
L|x-\hx|,\\
\ds |\Bf (t,0,u)| \leq L,
\end{array}\right.   \q  \all
(t,x,\hx,u)\in [0,T]\times \dbR^n\times \dbR^n\times U.
\end{equation}

Now, let $\bu(\cd)\in \cU_{ad}$  be a minimizer of $J(\cd)$ over
$\cU_{ad}$. We linearize $\cU_{ad}$ along $\bu(\cd)$ in the
following manner. Define
\begin{equation}\label{E205}
\cM_{ad}\equiv\set{ (1-\a)\d_{\bu(\cd)}+\a \d_{u(\cd)}\big| \a\in
[0,1], u(\cd)\in \cU_{ad}},
\end{equation}
where $\d_v$ denotes the Dirac measure at $v$ on $U$. For an element
$\s(\cd) = (1-\a)\d_{\bu(\cd)}+\a \d_{u(\cd)}\in \cM_{ad}$, denote
\begin{equation}\label{E206}
  \Bf(t,x,\s(t))\equiv \int_U \Bf(t,x,v)\s(t)(dv)= (1-\a)\Bf(t,x,\bu(t))+\a
\Bf(t,x, u(t)).
\end{equation}
Then we can define $x(\cd)=x(\cd;\s(\cd))$ as the solution of the
equation
\begin{equation}\label{E207}
\left\{\begin{array}{ll}\ds \dot x(t)=f(t,x(t),\s(t)),& \eqin [0,T],
\\
x(0)=x_0 &
\end{array}\right.
\end{equation}
and the corresponding cost functional $J(\s(\cd))$ by
\begin{equation}\label{E208}
J(\s(\cd))\equiv \int^T_0f^0(t,x(t;\s(\cd)),\s(t))\,  dt.
\end{equation}
We can see that $x(\cd;u(\cd))$ and $J(u(\cd))$ coincide  with
$x(\cd;\d_{u(\cd)})$ and
 $J(\d_{u(\cd)})$ respectively. Thus, $\cU_{ad}$ can be viewed as a subset of $\cM_{ad}$ in the sense of identifying
$u(\cd)\in \cU_{ad}$ to $\d_{u(\cd)}\in \cM_{ad}$. Readers who are
familiar with relaxed controls will
 immediately find $\cM_{ad}$ is a subset of relaxed control space. Yet,  elements
of $\cM_{ad}$ are much simpler than other relaxed controls. This is why we
need neither  to pose additional assumptions like that the control
domain is compact as Warga did (c.f. \cite{Wa3}) nor to introduce
the relaxed control defined by finite-additive probability measure
as Fattorini did (c.f. \cite{Fa1}).
 $\cM_{ad}$ has a linear structure at $\bu(\cd)$, i.e., it contains
 all elements in the form $\d_{\bu(\cd)}+\a (\d_{u(\cd)}-\d_{\bu(\cd)})$
($\a\in [0,1]$). It can be proved easily that $\d_{\bu(\cd)} $ is a
minimizer of $J(\s(\cd))$ over $\cM_{ad}$. Using this fact, we can
derive  Pontryagin's maximum principle from
\begin{equation}\label{E209}
0\leq \lim_{\a\to 0^+}{J((1-\a)\d_{\bu(\cd)}+ \a
 \d_{u(\cd)} )-J( \d_{\bu(\cd)})\over \a}.\end{equation}

It is easy to prove the following results.
\begin{Lemma}\tt\label{T201} Let \thb{S1}---\thb{S2} hold. Then, there exists a
constant $C>0$, such that for any  $\s(\cd)\in \cM_{ad}$,
\begin{equation}\label{E209B}
\left\{\begin{array}{l} \ds  \|x(\cd;\s(\cd))\|_{C[0,T]}\leq C, \\
 |x(t;\s(\cd))-x(\mt;\s(\cd))|\leq C|t-\mt|.\end{array}\right.
\end{equation}
\end{Lemma}
\Proof Let $\s(\cd)\in \cM_{ad}$ and $x(\cd)= x(t;\s(\cd))$. We have
$$
x(t)=x_0+\int^t_0f(s,x(s),\s(s))\, ds,\qq\all t\in [0,T].
$$
Then it follows from (S2) that
\begin{eqnarray*}\ds
\nnb |x(t)|&\leq & |x_0|+\ia{\int^t_0f(s,0,\s(s))\,
ds}+ L\int^t_0 |x(s)|\, ds\\
&\leq & \ds |x_0|+LT+L\int^t_0|x(s)|\, ds,\qq\all t\in [0,T].
\end{eqnarray*}
Thus, by Gronwall's inequality,
\begin{equation}\label{EX001}
|x(t)|\leq (|x_0|+LT)e^{Lt},\qq\all t\in [0,T].
\end{equation}
Consequently, by (S2),
$$
\ia{f(t,x(t),\s(t))}\leq L+L|x(t)|\leq L+L(|x_0|+LT)e^{LT},\qq\all
t\in [0,T].
$$
Therefore
\begin{equation}\label{EX002}\ds
|x(t)-x(\mt)|= \ia{\int^t_{\mt} f(s,x(s),\s(s))\, ds}\leq L \Big(1+
(|x_0|+LT)e^{LT}\Big) \ia{t-\mt}, \all t,\mt\in [0,T].
\end{equation}
We get \eqref{E209B} from \eqref{EX001}---\eqref{EX002}.
\endpf

\begin{Lemma}\tt\label{T202} Let \thb{S1}---\thb{S2} hold and
$\bu(\cd)$ be a minimizer of $J(u(\cd))$ over $\cU_{ad}$. Then
$\d_{\bu(\cd)} $ is a minimizer of $J(\s(\cd))$ over $\cM_{ad}$.
\end{Lemma}
\Proof Fix $\a\in [0,1]$ and $u(\cd)\in \cU_{ad}$. Denote
$\s^\a(\cd)=(1-\a) \d_{\bu(\cd)}+\a \d_{u(\cd)}$ and
$x^\a(\cd)=x(\cd;\s^\a(\cd))$. We will prove that
\begin{equation}\label{E210}
J( \d_{\bu(\cd)})\leq J(\s^\a(\cd)).
\end{equation}

For  $\ve>0$, define
$$
u^{\a,\ve}(t)=\left\{\begin{array}{ll}
\ds u(t),& \eqif  \set{{t\over \ve}} \in [0, \a), \\
\ds\bu(t),& \eqif  \set{{t\over \ve}} \in [\a,1),
\end{array}\right.
$$
 where $\{a\}$ denotes the decimal part of a real number $a$.
Then $u^{\a,\ve}(\cd)\in \cU_{ad}$. Denote
$x^{\a,\ve}(\cd)=x(\cd;u^{\a,\ve}(\cd))$. Then by Lemma \ref{T201},
$x^{\a,\ve}(\cd)$ is  uniformly bounded and equicontinuous on
$[0,T]$.  Consequently, by Arzel\'a-Ascoli's theorem, along a subsequence
$\ve\to 0^+$, $x^{\a,\ve}(\cd)$ converges uniformly to some $y(\cd)$ in $
 [0,T]$.  Thus,  using a generalization of Riemann-Lebesgue's Theorem
(see Ch. II, Theorem 4.15 in \cite{Zygmund}),  we can easily
prove that by a subsequence $\ve\to 0^+$,
\begin{equation}\label{EX003}
  f(t,y(t),u^{\a,\ve}(t))\to   f(t,y
(t),\s^\a (t)),  \q \mbox{weakly in }\, L^2(0,T;\dbR^n).
\end{equation}
Since
$$
\ia{f(t,x^{\a,\ve}(t),u^{\a,\ve}(t))-f(t,y(t),u^{\a,\ve}(t))}\leq
L\ia{x^{\a,\ve}(t)-y(t)},
$$
it follows from \eqref{EX003}   that
\begin{equation}\label{E211}
  f(t,x^{\a,\ve}(t),u^{\a,\ve}(t))\to   f(t,y
(t),\s^\a (t)),  \q \mbox{weakly in }\, L^2(0,T;\dbR^n).
\end{equation}
Similarly,
\begin{equation}\label{E212}
   f^0(t,x^{\a,\ve}(t),u^{\a,\ve}(t))\to   f^0(t,y
(t),\s^\a (t)), \q \mbox{weakly in }\, L^2(0,T).
\end{equation}
Passing to the limit for ~$\ve\to 0^+$ in the following equality
$$
x^{\a,\ve}(t)=x_0+\int^t_0 f(s,x^{\a,\ve}(s),u^{\a,\ve}(s))\, ds,
$$
we get from \eqref{E211}  that
$$
y(t)=x_0+\int^t_0 f(s,y(s),\s^\a (s))\, ds,
$$
i.e., $y(\cd)=x^\a(\cd)$. Furthermore, we can see that
  $x^{\a,\ve}(\cd)$ itself converges uniformly to
$x^\a(\cd)$ in $  [0,T]$. Combining this with \eqref{E212}, we have
\begin{equation}\label{E213}
J(\d_{\bu(\cd)})\leq \lim_{\ve\to 0^+} J(u^{\a,\ve}(\cd))=J(\s^\a
(\cd)).
\end{equation}
\endpf


\def\theequation{3.\arabic{equation}}
\setcounter{equation}{0} \setcounter{section}{3}
\setcounter{Definition}{0} \setcounter{Remark}{0}\textbf{3.
Pontryagin's Maximum Principle.} Now, we will derive Pontryagin's
maximum principle from \eqref{E209}. The idea  of our proof could be
tracked back to the works on relaxed control (c.f. \cite{Ga} and
\cite{Wa3}, for example). However, one can still find that the proof
we will give later has some improvement. Moreover, it can also be
used to problems governed by partial differential equations and even
having state constraints (c.f. \cite{Lou-Yong 2009}).

We keep the notations used in \S 2 and denote
$\bx(\cd)=x(\cd;\bu(\cd))$. We have
\begin{eqnarray}\label{E301}
\nnb   \ds  X^\a(t) &\equiv & {x^\a(t)-\bx(t)\over \a}   =
\int^t_0\Big[{f (s,x^\a
(s),\bu(s))-f (s,\bx(s),\bu(s))\over \a} \\
\nnb & & +
 f(s,x^\a(s),u(s))-f (s,x^\a(s),\bu(s)) \Big]\, ds\\
\nnb &=& \int^t_0\Big[\int^1_0  f_x(s,\bx(s)+\t (x^\a
(s)-\bx(s)),\bu(s))^\top   \, d\t \, X^\a(s) \\
 & & +
f(s,x^\a(s),u(s))-f (s,x^\a(s),\bu(s))  \Big]\, ds,
\end{eqnarray}
where  ~$f_x(t,x,u)$ denotes the transpose of the Jacobi matrix
of~$f$ on~$x$.
\if{
\begin{equation}\label{E302}
f_x(t,x,u)=\pmatrix{ {\pa f^1\over \pa  x_1}(t,x,u) & {\pa f^2\over
\pa  x_1}(t,x,u) & \ldots & {\pa f^n\over \pa  x_1}(t,x,u)\cr {\pa
f^1\over \pa  x_2}(t,x,u) & {\pa f^2\over \pa  x_2}(t,x,u) & \ldots
& {\pa f^n\over \pa  x_2}(t,x,u)\cr \vdots & \vdots & \ddots &
\vdots\cr {\pa f^1\over \pa  x_n}(t,x,u) & {\pa f^2\over \pa
x_n}(t,x,u) & \ldots & {\pa f^n\over \pa  x_n}(t,x,u)}.
\end{equation}
}\fi
By \eqref{E301}, (S2), and using the same argument as the proof of the uniform
convergence of $x^\a(\cd)\to \bx(\cd)$  in $[0, T]$, we can easily get
\begin{equation}\label{E303}
X^\a(\cd)\to X(\cd), \qq\mbox{uniformly in\, }  [0,T]
\end{equation}
and $X(\cd)$ solves the variational equation
\begin{equation}\label{E304}
\left\{\begin{array}{ll}\ds \dot X(t)=  f_x(t,\bx(t),\bu(t))^\top \,
X(t) + f(t,\bx (t),u(t))-f (t,\bx(t),\bu(t)) ,& \eqin [0,T],
\\
X(0)= 0. &
\end{array}\right.
\end{equation}
Now, by introducing the adjoint equation
\begin{equation}\label{E305}
\left\{\begin{array}{ll}\ds \dot \bpsi (t)=- f_x(t,\bx(t),\bu(t)) \,
\bpsi(t) + f^0_x(t,\bx (t),\bu(t))  ,& \eqin [0,T],
\\
\bpsi(T)= 0, &
\end{array}\right.
\end{equation}
we get from \eqref{E213} and Lebesgue's dominated convergence
theorem that
\begin{eqnarray}\label{E306}
\nnb 0 &  \leq & \ds \lim_{\a\to 0^+} {J((1-\a)\d_{\bu(\cd)}+ \a
 \d_{u(\cd)} )-J( \d_{\bu(\cd)})\over \a}\\
\nnb &=& \lim_{\a\to 0^+} \int^T_0\Big[{f^0(t,x^\a
(t),\bu(t))-f^0(t,\bx(t),\bu(t))\over \a} \\
\nnb & & \q +
 f^0(t,x^\a(t),u(t))-f^0(t,x^\a(t),\bu(t)) \Big]\, dt\\
\nnb &=& \lim_{\a\to 0^+} \int^T_0\Big[\int^1_0
\ip{f^0_x(t,\bx(t)+s(x^\a
(t)-\bx(t)),\bu(t)), X^\a(t)}\, ds \\
\nnb & & \q +
 f^0(t,x^\a(t),u(t))-f^0(t,x^\a(t),\bu(t)) \Big]\, dt\\
\nnb &= & \int^T_0 \Big[\ip{f^0_x(t,\bx(t),\bu(t)), X(t)}+
 f^0(t,\bx (t),u(t))-f^0(t,\bx (t),\bu(t)) \Big]\, dt\\
  &=& \int^T_0 \Big[H(t,\bx(t),\bu(t)),\bpsi(t))-
H(t,\bx(t),u(t)),\bpsi(t))\Big]\, dt,
\end{eqnarray}
where
\begin{equation}\label{E307}
H(t,x,u,\psi)\equiv \ip{f(t,x,u),\psi}-f^0(t,x,u),\qq\all
(t,x,u,\psi)\in [0,T]\times \dbR^n\times U\times \dbR^n.
\end{equation}
Then, since $U$ is separable and $H$ is continuous in $u$, it
follows from (\ref{E306}) and a  standard argument that
\begin{equation}\label{E308}
H(t,\bx(t),\bu(t)),\bpsi(t))=\max_{v\in U} H(t,\bx(t),v,\bpsi(t)),
\qq\eqae t\in [0,T].
\end{equation}
 Relations  \eqref{E305}, \eqref{E307}, \eqref{E308} form  Pontryagin's maximum
principle.

\bigskip

\def\theequation{4.\arabic{equation}}
\setcounter{equation}{0} \setcounter{section}{4}
\setcounter{Definition}{0} \setcounter{Remark}{0}\textbf{4.
Second-Order Necessary Optimality Conditions.} We turns to study
second-order necessary optimality conditions where Pontryagin's maximum principle is viewed as
a first-order necessary optimality
condition. In other words, we will give a second-order necessary
condition for optimality to distinguish singular controls in the
sense of Pontryagin's maximum principle. One can see that in
\eqref{E308}, the equality holds if and only if
\begin{equation}\label{E401}
u(t)\in \bU(t)\equiv \set{w\big| H(t,\bx(t),w,\bpsi(t))=\max_{v\in
U} H(t,\bx(t),v,\bpsi(t))}, \qq\eqae t\in [0,T].
\end{equation}
In this case,
\begin{equation}\label{E402}
\int^T_0 \Big[\ip{f^0_x(t,\bx(t),\bu(t)), X(t)}+
 f^0(t,\bx (t),u(t))-f^0(t,\bx (t),\bu(t)) \Big]\, dt=0.
\end{equation}
Denote
\begin{equation}\label{E402A}
\bU_{ad}= \set { v(\cd)\in \cU_{ad}| v(t) \in \bU(t), \q\eqae t\in
[0,T]\,}.
\end{equation}
Elements in $\bU_{ad}$ are called  singular controls in the sense of
Pontryagin's maximum principle. If $U$ is an open subset of
$\dbR^m$, then
\begin{equation}\label{ES001}
H_u(t,\bx(t),\bu(t),\bpsi(t))=0, \qq\eqae t\in [0,T].
\end{equation}
In this case, we call elements in
\begin{equation}\label{ES002}
\set{v(\cd)\in \cU_{ad}\big| H_u(t,\bx(t),v(t),\bpsi(t))=0, \q \eqae t\in [0,T]}
\end{equation}
as singular controls in the classical sense (see Definitions 1 and 2
in \cite{G-K}).

Now we  make the following assumption:

(S3) Functions $\Bf$ are twice continuously differentiable in $x$.
Moreover,  it holds that
\begin{equation}\label{E402B}
|\Bf_x(t,x,u)-\Bf_x (t,\hx, u)|\leq L|x-\hx|, 
\q  \all (t,x,\hx,u)\in [0,T]\times \dbR^n\times \dbR^n\times U.
\end{equation}
We mention that (S2) implies
$$
|\Bf_x(t,x,u)|\leq L, \q  \all (t,x, u)\in [0,T] \times \dbR^n\times
U
$$
and (S3) implies
$$
|f^k_{xx}(t,x,u)|\leq L, \q  \all (t,x, u)\in [0,T] \times
\dbR^n\times U
$$
for $k=0,1,2,\ldots,n$.

 The following theorem gives second-order
necessary optimality conditions.
\begin{Theorem}\tt\label{T401}  Let \thb{S1}---\thb{S3} hold and $\bu(\cd)$ be a
minimizer of $J(\cd)$ over $\cU_{ad}$. Define $\bW(\cd)$ be the
solution of the following second-order adjoint equation:
\begin{equation}\label{E409}
\left\{\begin{array}{ll}\ds \dot  \bW  (t)+f_x(t,\bx(t),\bu(t))
\bW(t)+\bW(t) f_x(t,\bx(t),\bu(t))^\top\\
\qq\qq\qq\qq\ds +H_{xx}(t,\bx(t),\bu(t),\bpsi(t))=0, \q \eqin [0,T],\\
\bW(T)=0 \end{array}\right.
\end{equation}
and $\bGP(\cd)$ be the solution of
\begin{equation}\label{E409B}
\left\{\begin{array}{ll}\ds \dot \bGP(t)=  f_x(t,\bx(t),\bu(t))^\top
\, \bGP(t),& \eqin [0,T],
\\
\bGP(0)= I, &
\end{array}\right.
\end{equation}
where $I$ is the unit $n\times n$ matrix. Then for any $u(\cd)\in
\bU_{ad}$,
\begin{eqnarray}\label{E411A}
\nnb  &    & \ds \int^T_0dt \int^t_0\Big\langle
\bW(t)(f(t,\bx(t),\bu(t))-f(t,\bx(t),
 u(t)))\\
 \nnb & & \qq\q + H_x(t,\bx(t),\bu(t),\bpsi(t))-H_x(t,\bx(t),
u(t),\bpsi(t)),\\
&&\qq \q\bGP(t)\bGP(s)^{-1}(f(s,\bx(s),\bu(s))-f(s,\bx(s),
 u(s)))\Big\rangle\, ds\leq 0.
\end{eqnarray}
\end{Theorem}
\Proof Let $u(\cd)\in \bU_{ad}$. Then by \eqref{E402},  for any
$\a\in (0,1]$,
\begin{eqnarray}\label{E403}
\nnb  &   & \ds   J((1-\a)\d_{\bu(\cd)}+ \a
 \d_{u(\cd)} )-J( \d_{\bu(\cd)}) \\
\nnb &=&   \a \int^T_0\Big[\int^1_0 \ip{f^0_x(t,\bx(t)+s(x^\a
(t)-\bx(t)),\bu(t)), X^\a(t)}\, ds \\
\nnb & &  \qq +
 f^0(t,x^\a(t),u(t))-f^0(t,x^\a(t),\bu(t)) \Big]\, dt\\ \if{
\nnb &=&  \a  \int^T_0 dt \int^1_0  \ip{f^0_x(t,\bx(t)+s(x^\a
(t)-\bx(t)),\bu(t)), X^\a(t)}-\ip{f^0_x(t,\bx(t),\bu(t)), X(t)}\Big]\, ds \\
\nnb & &  +\a\int^T_0  \big(
 f^0(t,x^\a(t),u(t))-f^0(t,\bx (t),u(t))\big)-\big(
 f^0(t,x^\a(t),\bu(t))-f^0(t,\bx (t),\bu(t)) \big) \Big]\, dt\\}\fi
\nnb &=&   \a \int^T_0dt\int^1_0  \ip{f^0_x(t,\bx(t)+s(x^\a
(t)-\bx(t)),\bu(t))-f^0_x(t,\bx(t),\bu(t)), X^\a(t)} \, ds \\
\nnb && + \a\int^T_0    \ip{f^0_x(t,\bx(t),\bu(t)), X^\a(t) -
X(t)} \, dt  \\
\nnb & &  +\a \int^T_0\Big[\big(
 f^0(t,x^\a(t),u(t))-f^0(t,\bx (t),u(t))\big)-\big(
 f^0(t,x^\a(t),\bu(t))-f^0(t,\bx (t),\bu(t)) \big) \Big]\, dt\\
\nnb &=&   \a^2 \int^T_0dt\int^1_0 ds\int^1_0
 s \ip{f^0_{xx} (t,\bx(t)+s\t (x^\a
(t)-\bx(t)),\bu(t))X^\a(t), X^\a(t)} \, d\t \\
\nnb && + \a^2\int^T_0     \ip{f^0_x(t,\bx(t),\bu(t)),
Y^\a(t)} \, dt \\
\nnb & &  +\a^2 \int^T_0 dt\int^1_0  \ip{
 f^0_x(t,\bx(t)+s  (x^\a
(t)-\bx(t)),u(t)), X^\a(t)}\, ds \\
  & &  -\a^2 \int^T_0 dt\int^1_0 \ip{
 f^0_x(t,\bx(t)+s  (x^\a
(t)-\bx(t)),\bu(t)),X^\a(t)}  \, ds,
\end{eqnarray}
where $\ds Y^\a(\cd)\equiv \pmatrix{Y^\a_1(\cd) & Y^\a_2(\cd) &
\ldots & Y^\a_n(\cd)}^\top={X^\a(\cd)- X(\cd)\over \a}$ satisfies
\begin{eqnarray}\label{E404}
\nnb   Y^\a_k(t) & = & \ds {1\over \a}\int^t_0\Big[\int^1_0
\ip{f^k_x(s,\bx(s)+\t (x^\a
(s)-\bx(s)),\bu(s)), X^\a(s)} \, d\t \,  \\
\nnb  & & +
f^k(s,x^\a(s),u(s))-f^k (s,x^\a(s),\bu(s))  \Big]\, ds\\
\nnb & & -{1\over \a}\int^t_0\Big[ \ip{f^k_x(s,\bx(s),\bu(s)), X
(s)}   +
f^k(s,\bx ,u(s))-f^k (s,\bx (s),\bu(s))  \Big]\, ds\\
\nnb &=&   \int^t_0ds\int^1_0 d\t\int^1_0
 \t \ip{f^k_{xx} (s,\bx(s)+\t\z (x^\a
(s)-\bx(s)),\bu(s))X^\a(s), X^\a(s)} \, d\z \\
\nnb & &  +  \int^t_0 ds\int^1_0  \ip{
 f^k_x(s,\bx(s)+\t  (x^\a
(s)-\bx(s)),u(s)), X^\a(s)}\, d\t \\
\nnb  & &  -  \int^t_0 ds\int^1_0 \ip{
 f^k_x(s,\bx(s)+\t  (x^\a
(s)-\bx(s)),\bu(s)),X^\a(s)}  \, d\t\\
  && +  \int^t_0 \ip{f^k_x(s,\bx(s),\bu(s)), Y^\a(s)} \, ds ,\q
k=1,2,\ldots,n.
\end{eqnarray}
Using (S3) and by the same way to derive
 \eqref{E303}, we can
get
\begin{equation}\label{E405} Y^\a(\cd)\to Y(\cd),
\qq\mbox{uniformly in\, }  [0,T]
\end{equation}
where $Y(\cd)$ solves the following second variational equation
\begin{equation}\label{E406}
\left\{\begin{array}{l} \ds \dot Y(t) =  f_x(t,\bx(t),\bu(t))^\top Y
(t)+(f
_x(t,\bx(t),u(t))-f _x(t,\bx(t), \bu(t)))^\top X(t)\\
\\
  \qq\qq\qq\qq\qq  +
 {1\over 2} \pmatrix{ \ip{f^1_{xx} (t,\bx(t),\bu(t))X (t), X (t)} \cr
\ip{f^2_{xx} (t,\bx(t),\bu(t))X (t), X (t)}\cr  \vdots\cr
\ip{f^n_{xx} (t,\bx(t),\bu(t))X (t), X (t)}}, \q \eqin [0,T],\\
Y(0) =0.    \end{array}\right.
\end{equation}
Then it follows from \eqref{E403} and Lebesgue's dominated
convergence theorem  that
\begin{eqnarray}\label{E407}
\nnb 0 &\leq    & \ds   \lim_{\a\to 0^+}{J((1-\a)\d_{\bu(\cd)}+ \a
 \d_{u(\cd)} )-J( \d_{\bu(\cd)}) \over \a^2}\\
 \nnb &=&
\ \int^T_0     \ip{f^0_x(t,\bx(t),\bu(t)), Y (t)} \, dt    +
\int^T_0   \ip{
 f^0_x(t,\bx(t),u(t))- f^0_x(t,\bx(t) ,\bu(t)), X (t)}\, dt\\
\nnb & & +{1\over 2}  \int^T_0
  \ip{f^0_{xx} (t,\bx(t),\bu(t))X (t), X (t)} \, d t\\
\nnb &=& \int^T_0\ip{H_x(t,\bx(t),\bu(t),\bpsi(t))-H_x(t,\bx(t),
u(t),\bpsi(t)),X(t)}\, dt\\
\nnb & & - {1\over
2}\int^T_0\ip{H_{xx}(t,\bx(t),\bu(t),\bpsi(t))X(t),X(t)}\, dt\\
\nnb &=& \int^T_0\ip{H_x(t,\bx(t),\bu(t),\bpsi(t))-H_x(t,\bx(t),
u(t),\bpsi(t)),X(t)}\, dt\\
  & & - {1\over 2}\tr \Big[\int^T_0
H_{xx}(t,\bx(t),\bu(t),\bpsi(t))X(t) X(t)^\top \, dt\Big],
\end{eqnarray}
where $\tr B$ denotes the trace of a matrix $B$.
One can easily verify that
\begin{equation}\label{E408}
\left\{\begin{array}{l} \ds  \doo {} t\big(X(t)X(t)^\top\big)  =
f_x(t,\bx(t),\bu(t))^\top
\big(X(t)X(t)^\top\big)+\big(X(t)X(t)^\top\big)f_x(t,\bx(t),\bu(t))\\
\qq\qq\qq\qq\ds +(f(t,\bx(t),u(t))-f(t,\bx(t), \bu(t)))X(t)^\top\\
\qq\qq\qq\qq\ds+X(t)(f(t,\bx(t),u(t))-f(t,\bx(t), \bu(t)))^\top, \q \eqin [0,T],\\
X(0)X(0)^\top =0.   \end{array}\right.
\end{equation}
Now, we introduce the second-order adjoint equation \eqref{E409}. By (S1)---(S3), we can see that
\eqref{E409} admits a unique solution $\bW(\cd)$. Since  $\bW(\cd)^\top$ satisfies \eqref{E409} too,
 $\bW(\cd)$ should be symmetric. Since  $\ds \tr
(AB)=\tr(BA)$ for all $k\times j$ matrix $A$ and $j\times k$ matrix $B$, we have
\begin{eqnarray}\label{E410}
\nnb   & &\ds  - {1\over 2}\tr \Big[\int^T_0
H_{xx}(t,\bx(t),\bu(t),\bpsi(t))X(t) X(t)^\top \, dt\Big]\\
\nnb &=&    {1\over 2}\tr \Big[\int^T_0
\Big(\doo {} t\bW (t)+f_x(t,\bx(t),\bu(t)) \bW(t)+\bW(t)
f_x(t,\bx(t),\bu(t))^\top\Big)X(t) X(t)^\top \, dt\Big]\\
\nnb &=&   {1\over 2}\tr \Big\{\int^T_0 \Big[-\bW(t)\doo {} t\big(X(t)
X(t)^\top\big)  +f_x(t,\bx(t),\bu(t)) \bW(t) X(t)
X(t)^\top \\
\nnb & & +\bW(t) f_x(t,\bx(t),\bu(t))^\top X(t)
X(t)^\top \Big]\, dt\Big\}\\
\nnb &=&    {1\over 2}\tr \Big\{\int^T_0
\Big[-\bW(t)f_x(t,\bx(t),\bu(t))^\top
 X(t)X(t)^\top -\bW(t) X(t)X(t)^\top f_x(t,\bx(t),\bu(t))\\
\nnb & &  -\bW(t)(f(t,\bx(t),u(t))-f(t,\bx(t), \bu(t)))X(t)^\top\\
\nnb & &  -\bW(t)X(t)(f(t,\bx(t),u(t))-f(t,\bx(t), \bu(t)))^\top\\
\nnb & &  +f_x(t,\bx(t),\bu(t)) \bW(t) X(t) X(t)^\top +\bW(t)
f_x(t,\bx(t),\bu(t))^\top X(t) X(t)^\top \Big]\, dt\Big\}\\
 &=&  \int^T_0 \ip{ \bW(t)(f(t,\bx(t),\bu(t))-f(t,\bx(t),
 u(t))),X(t)}\, dt.
\end{eqnarray}
By \refeq{E304} and \refeq{E409},
$$
X(t)=\int^t_0\bGP(t)\bGP(s)^{-1}
 (f(s,\bx(s),u(s))-f(s,\bx(s),
 \bu(s))) \, ds.
$$
Thus, it follows from \eqref{E407} and \refeq{E410} that
\begin{eqnarray*}
\nnb 0 &\leq    & \ds \int^T_0\Big\langle
\bW(t)(f(t,\bx(t),\bu(t))-f(t,\bx(t),
 u(t)))\\
 \nnb & &\ds + H_x(t,\bx(t),\bu(t),\bpsi(t))-H_x(t,\bx(t),
u(t),\bpsi(t)),X(t)\Big\rangle\, dt\\
\nnb &=& \ds \int^T_0dt \int^t_0\Big\langle
\bW(t)(f(t,\bx(t),\bu(t))-f(t,\bx(t),
 u(t)))+ H_x(t,\bx(t),\bu(t),\bpsi(t))\\
 \nnb & & \ds -H_x(t,\bx(t),
u(t),\bpsi(t)), \bGP(t)\bGP(s)^{-1}
 (f(s,\bx(s),u(s))-f(s,\bx(s),
 \bu(s)))\Big\rangle\, ds.
\end{eqnarray*}
Therefore, we  finish the proof.
\endpf
\if{
\begin{Remark}\tt
If we suppose that $U$ is a domain in $\dbR^m$ and $f^0,f$ are twice continuously differentiable
in $x$ and $u$, then we can get from Theorem~ \ref{T401} that
\begin{equation}\label{E002}
 \int^T_0\Big[\ip{
H_{xx}(t,\bx(t),\bu(t),\bpsi(t))z(t),z(t)}+2\ip{H_{xu}(t,\bx(t),\bu(t),\bpsi(t))v(t),z(t)}\Big]\, dt\leq 0,
\end{equation}
\if{
 (set $u=\bu+\b v$)
\begin{eqnarray*}
0 &\leq & \lim_{\b \to 0} {1\over \b^2}\Big[\int^T_0\Big\langle
\bW(t)(f(t,\bx(t),\bu(t))-f(t,\bx(t),
 u(t)))\\
 \nnb & &\ds\qq + H_x(t,\bx(t),\bu(t),\bpsi(t))-H_x(t,\bx(t),
u(t),\bpsi(t)),X(t)\Big\rangle\, dt\Big]\\
&=& -\int^T_0\Big\langle
\bW(t) f_u(t,\bx(t),\bu(t))v(t) + H_{xu}(t,\bx(t),\bu(t),\bpsi(t))v(t),z(t)\Big\rangle\, dt\\
&=& -{1\over 2}\int^T_0\Big[\ip{
H_{xx}(t,\bx(t),\bu(t),\bpsi(t))z(t),z(t)}+2\ip{H_{xu}(t,\bx(t),\bu(t),\bpsi(t))v(t),z(t)}\Big]\, dt,
\end{eqnarray*}}\fi
where
$$
\left\{\begin{array}{ll}\ds \dot z(t)=f_x(t,\bx(t),\bu(t))z(t)+f_u(t,\bx(t),\bu(t))v(t),& \eqin [0,T],
\\
z(0)=0. &
\end{array}\right.
$$
On the other hand, the
classical results appear as (see \cite{Ros}, for example)
\begin{equation}\label{E003}
\int^T_0 \Big[\ip{H_{xx} z, z}+2\ip{H_{xu}v, z}+\ip{H_{uu}v,v}\Big]\, dt\leq 0.
\end{equation}
Since it follows from the first-order necessary condition \refeq{E308} that ~$$H_{uu}(t,\bx(t),\bu(t),\bpsi(t))\leq 0,\qq\eqae t\in [0,T],$$
we can see that \refeq{E002} implies \refeq{E003}.
\end{Remark}
}\fi

To get an analogue of  the maximum condition
in Pontryagin's maximum principle, we introduce a lemma concerned
with the well-known Filippov's Lemma or measurable selection. We
recall that a Polish space is a  separable completely metrizable topological space. We mention that
all (nonempty) closed sets and open sets in $\dbR^m$ are polish
spaces.
\begin{Lemma}\tt\label{T402} Let $X$ be a Polish space,
$\cT\subseteq \dbR^n$ be a Lebesgue measurable set. Assume that
$\G:\cT\to 2^X$ is measurable  \thb{i.e., for any closed set $F$,
$\set{t\in \cT|\G(t)\in F}$ is measurable} and takes values on the
family of nonempty closed subsets of ~$X$. Then $\G(\cd)$ admits a measurable selection,
i.e., there exists a Lebesgue measurable map $\g: \cT\to X$, such
that
$$\g(t)\in \G(t), \q \eqae t\in \G(t).$$
\end{Lemma}
Lemma \ref{T402} and its proof  can be found in \cite{Li-Yong}
(Ch. 3, Theorem 2.23).  Based on this lemma, we have that
 \begin{Theorem}\tt\label{T403} Let  $U$ be a Polish space. Then under assumptions of Theorem
\ref{T401}, for almost all $t\in [0,T]$, it holds that:
\begin{eqnarray}\label{E414A}
\nnb  &    & \ds  \Big\langle \bW(t)(f(t,\bx(t),\bu(t))-f(t,\bx(t),
 v))+ H_x(t,\bx(t),\bu(t),\bpsi(t))\\
   & & -H_x(t,\bx(t),
v,\bpsi(t)), (f(t,\bx(t),\bu(t))-f(t,\bx(t),
 v))\Big\rangle \leq 0, \q \all v\in \bU(t).
\end{eqnarray}
\end{Theorem}
\Proof  Denote
\begin{eqnarray*}
\nnb F(t,u) &  =  & \ds  \bGP(t)^\top \Big[
 \bW(t)(f(t,\bx(t),\bu(t))-f(t,\bx(t),
 u )) \\
&& +  H_x(t,\bx(t),\bu(t),\bpsi(t)) -H_x(t,\bx(t),
u ,\bpsi(t))\Big],\\
G(t,u) &  =  & \ds \bGP(t)^{-1}
 (f(t,\bx(t),u )-f(t,\bx(t),
 \bu(t))) ,\\
 && \qq\qq\qq\qq (t,u)\in [0,T]\times U.
\end{eqnarray*}
Then by Theorem \ref{T401}, for any $u(\cd)\in \bU_{ad}$, we have
\begin{eqnarray*}
\nnb  &    & \ds \int^T_0dt \int^t_0 \ip{F(t,u(t)),G(s,u(s))}\,
ds\leq 0.
\end{eqnarray*}
Let
\begin{equation}\label{E414}
  E_{u(\cd)}=\set{ t\in [0,T)\big| \lim_{\a\to 0^+}{1\over
 \a} \int^{t+\a}_t \pmatrix{F(s,u(s))\cr |F(s,u(s))|^2\cr G(s,u(s))}\, ds=\pmatrix{F(t,u(t))\cr |F(t,u(t))|^2\cr G(t,u(t))} }.
\end{equation}
Then $E_{u(\cd)} $  has Lebesgue measure $T$. Let $\b\in
E_{u(\cd)}$. For $\a\in (0, T-\b) $, define
$$
u^\a(t)=\left\{\begin{array}{ll}\ds  \bu(t), & \eqif t\not\in
[\b ,\b+\a],\\
\ds u(t), & \eqif t\in [\b ,\b+\a].\end{array}\right.
$$
Then, $u^\a(\cd) \in \bU_{ad}$ and
\begin{equation}\label{E415}
  \int^{\b+\a}_\b  dt \int^t_\b \ip{ F(t,u(t)),G(s,u(s))}\, ds=\int^T_0dt \int^t_0 \ip{F(t,u^\a(t)),G(s,u^\a(s))}\,
ds\leq 0.
\end{equation}
It is easy to see that under assumptions (S1)---(S3),
$F(\cd,\cd),G(\cd,\cd)$ are uniformly bounded. Thus
$$
|F(t,u)|+|G(t,u)|\leq C, \qq \all t\in [0,T]\times U
$$
for some constant $C>0$. Consequently,
\begin{eqnarray*}
\nnb   &     & \ds  \Big| {2\over \a^2}\int^{\b+\a}_\b
dt \int^t_\b \ip{F(t,u(t)),G(s,u(s))}\, ds-\ip{F(\b,u(\b)),G(\b,u(\b))}\Big|\\
\nnb &\leq & \Big| {2\over \a^2}\int^{\b+\a}_\b \ip{
(t-\b)F(t,u(t)), \,
 {1\over t-\b}\int^t_\b G(s,u(s))\, ds-G(\b,u(\b)) }\, dt\Big|  \\
\nnb & & +\Big|  \ip{ {2\over \a^2}\int^{\b+\a}_\b
 (t-\b)(F(t,u(t))-F(\b,u(\b))) \, dt, \, G(\b,u(\b))} \Big|\\
\nnb &\leq &   {2C\over \a^2}\int^{\b+\a}_\b   (t-\b)\, dt  \,
 \sup_{r\in (\b,\b+\a]}\Big|{1\over r-\b}\int^r_\b G(s,u(s))\, ds-G(\b,u(\b)) \Big|  \\
\nnb & & +   2C\Big[{1\over \a^3}\int^{\b+\a}_\b
 (t-\b)^2\, dt\Big]^{1\over 2} \, \Big[{1\over \a}\int^{\b+\a}_\b\Big| F(t,u(t))-F(\b,u(\b)) \Big|^2\, dt \Big]^{1\over 2}  \\
\nnb &=&   C
 \sup_{r\in [\b,\b+\a]}\Big|{1\over r-\b}\int^r_\b G(s,u(s))\, ds-G(\b,u(\b)) \Big|  \\
\nnb & & +   {2\sqrt 3 C\over 3} \, \Big\{{1\over \a}\int^{\b+\a}_\b\Big[ |F(t,u(t))|^2+ | F(\b,u(\b))|^2-2\ip{F(t,u(t)),F(\b,u(\b))} \Big]\, dt \Big\}^{1\over 2}  \\
\end{eqnarray*}
Let $\a\to 0^+$, we get from \eqref{E414} that
$$
\lim_{\a\to 0^+}{2\over \a^2}\int^{\b+\a}_\b dt \int^t_\b
\ip{F(t,u(t)),G(s,u(s))}\, ds=\ip{F(\b,u(\b)),G(\b,u(\b))}.
$$
Therefore, combining the above with \eqref{E415}, we get
\begin{equation}\label{E416}
     \ip{F(\b,u(\b)),G(\b,u(\b))}\leq 0.
\end{equation}
That is, for any $u(\cd)\in \bU_{ad}$,
\begin{equation}\label{E417}
\nnb    \ip{F(t,u(t)),G(t,u(t))}\leq 0, \qq\eqae t\in [0,T].
\end{equation}
For $k=1,2,\ldots$, denote
$$
\cT_k=\set{t\in [0,T]\Big|\, \exists v\in \bU(t), \eqst
\ip{F(t,v),G(t,v)}\geq {1\over k}}.
$$
Then, $\cT_k$ is measurable. We claim $\cT_k$ has zero measure.
Otherwise, for any $t\in \cT_k$, $$\ds \G_k(t)\equiv \set{v\in
\bU\big|\ip{F(t,v),G(t,v)}\geq {1\over k}}$$ is a nonempty closed
subset of $U$. It is easy to see that $\G(\cd)$ is measurable
since $\ds \ip{F(t,v),G(t,v)}$ is measurable in $t$ and continuous
in $v$. Thus, by Lemma \ref{T402}, there  exists a measurable
function $\tiu_k:\cT_k\to \G_k$. Define
$$
u_k(t)=\left\{\begin{array}{ll}\ds \bu(t), & \eqif t\not\in \cT_k,\\
\ds \tiu_k(t), & \eqif t\in \cT_k.\end{array}.\right.
$$
Then $u_k(\cd)\in \bU_{ad}$ and
$$
\ip{F(t,u_k(t)),G(t,u_k(t))}\geq {1\over k}, \qq\eqae t\in \cT_k.
$$
Contradict to \eqref{E417}. Therefore $\cT_k$ has zero measure.
Consequently,
$$
\cT\equiv \set{t\in [0,T]\Big|\, \exists v\in \bU(t), \eqst
\ip{F(t,v),G(t,v)}>0}=\bigcup_{k=1}^\infty \cT_k
$$
has zero measure too. That is, for almost all $t\in [0,T]$,
\begin{equation}\label{E420}
\ip{F(t,v),G(t,v)}\leq 0, \qq\all v\in \bU(t).
\end{equation}
This completes the proof.
\endpf

It is not necessary to suppose that $U$ is a Polish space
in yielding \eqref{E417}. However, usually, we can not get
\eqref{E420} from \eqref{E417} if we only suppose that
 $U$ is a separable metric space. To see this, we introduce the
following example.
\begin{Example}\tt Let $T>0$ and $U\subset [0,T]$ be a non-measurable
set which contains $0$ and has no any subset of positive
measure. Let
$$
H(t,u)=-u^2(t-u)^2, \q g(t,u)=u^2, \qq (t,u)\in [0,T]\times U.
$$
Then $H$ and $g$ are smooth. Define
$$
\bU(t)\equiv \set{u\in U| H(t,u)=\max_{v\in U} H(t,v)},\qq t\in
[0,T].
$$
Then
$$
\bU(t)=\left\{\begin{array}{ll}\ds  \set{0} & \eqif t\not\in U,\\
\ds \set{0,t}, &\eqif t\in U.
\end{array}\right.
$$
Thus, if $$u(\cd)\in \bU_{ad}\equiv \set{v:[0,T]\to U|v(\cd)\,
\mbox{measurable}, v(t)\in \bU(t), \eqae [0,T]},$$  we must have
$u(t)=0, \eqae t\in [0,T]$. Otherwise, if $E\equiv \set{t|u(t)\ne
0}$ has positive measure, then $E\subseteq U$. This contradicts to
the assumption that U has no any subset of positive measure

Therefore,  for any $u(\cd)\in \bU_{ad}$,
\begin{equation}\label{E421}
g(t,u(t))\leq 0, \qq \eqae t\in [0,T].
\end{equation}
However,
$$
\set{t\in [0,T]\big|\max_{v\in \bU(t)} g(t,v)>0}=U\setminus \set{0}
$$
has not  zero measure. This means that the following statement does
not hold: for almost all $t\in [0,T]$,
$$
g(t,v)\leq 0,\qq \all v\in \bU(t).
$$
\end{Example}

\bigskip

\def\theequation{5.\arabic{equation}}
\setcounter{equation}{0} \setcounter{section}{5}
\setcounter{Definition}{0} \setcounter{Remark}{0}\textbf{5.
Sufficient Conditions.} Now, we give a sufficient condition for a
control being a local minimizer. In addition to (S1)---(S3), we
suppose that

(S4) There exists a modulus of continuity $\o: [0,+\infty)\to
[0,+\infty)$, such that
\begin{equation}\label{EX501}
\left\{\begin{array}{l} \ds |\Bf (t,x,u)-\Bf  (t, x,\hu)|\leq
 \o \big(\rho(u,\hu)\big),\\
\ds |\Bf_x (t,x,u)-\Bf_x (t, x,\hu)|\leq
 \o \big(\rho(u,\hu)\big),
 \end{array}\right.   \q  \all
(t,x,\hx,u)\in [0,T]\times \dbR^n\times \dbR^n\times U
\end{equation}
and  for
 $k=0,1,2,\ldots, n$, it holds that
\begin{equation}\label{EX502}
\begin{array}{l}\ds
  \ds |f^k_{xx}(t,x,u)-f^k_{xx}
(t,\hx,\hu)|\leq \o\big (|x-\hx|+\rho(u,\hu)\big), \\
\ds \qq\qq\qq\qq \q  \all (t,x,\hx,u,\hu)\in [0,T]\times
\dbR^n\times \dbR^n\times U\times U.\end{array}
\end{equation}

We note that (S2)---(S4) imply
$$
\left\{\begin{array}{l} \ds |\Bf (t,x,u)-\Bf  (t, \hx,\hu)|\leq
L|x-\hx|+
 \o \big(\rho(u,\hu)\big),\\
\ds |\Bf_x (t,x,u)-\Bf_x (t, \hx,\hu)|\leq L|x-\hx|+
 \o \big(\rho(u,\hu)\big).
 \end{array}\right.
$$

For $\bu(\cd),u(\cd)\in \cU_{ad}$ and $\a\in (0,1]$, we keep the
notations used in Sections 3 and 4. For notation simplicity,
denote
\begin{equation}\label{EX503}
\Th(t)=\o\big(\rho(u(t),\bu(t))\big), \qq t\in [0,T]
\end{equation}
and  $x^\a(\cd)=\bx(\cd), X^\a(\cd)=X(\cd)$, $ Y^\a(\cd)=Y(\cd)$
when $\a=0$.

The following lemma gives estimates of $X^\a(\cd)$ and $Y^\a(\cd)$.
\begin{Lemma}\tt\label{T501}  Assume \thb{S1}---\thb{S4},  $\bu(\cd),u(\cd)\in \cU_{ad}$ and $\a\in
[0,1]$. Then for any $\a\in [0,1]$,
\begin{equation}\label{E501}
|X^\a(t)|\leq C\int^t_0 \Th(s) \, ds, \qq\all x\in [0,T],
\end{equation}
\begin{equation}\label{E502}
|X^\a(t)|^2\leq C\int^t_0 \Big[\Th(s)\Big]^2\, ds, \qq\all x\in
[0,T],
\end{equation}
\begin{equation}\label{E503}
  |Y^\a(t)|  \leq   C  \int^t_0 \Big[\Th(s)\Big]^2\,
ds,\qq\all t\in [0,T]
\end{equation}
and
\begin{eqnarray}\label{E504}
 \nnb & &  |Y^\a(t)-Y(t)|  \\
 & \leq  & C\Big[\o\Big(C\int^T_0 \Th(s)\, ds \Big)   +   \int^T_0 \Th(s)\, ds \Big]\,
\int^t_0 \Big[\Th(s)\Big]^2\, ds\,,\qq\all t\in [0,T],
\end{eqnarray}
where and hereafter, $C>0$ denotes a constant, which is independent
of $u(\cd)$, $\a\in [0,1]$, and may be different in different lines.
\end{Lemma}
\Proof Let  $\a\in [0,1]$. By (S2)---(S4),  \eqref{E301}  and
\eqref{E304},
\begin{equation}\label{E505}
|X^\a(t)|\leq \int^t_0 \Big[ L|X^\a(s)|+\Th(s)\Big]\, ds,\qq\all
t\in [0,T].
\end{equation}
Then, \eqref{E501} follows from Gronwall's  inequality.
While
\eqref{E502} follows from \eqref{E501} and Cauchy-Schwarz's
inequality.

Similarly, by (S2)---(S4), \eqref{E404} and \eqref{E406},
\begin{equation}\label{E506}
|Y^\a(t)|\leq C\int^t_0 \Big[ |Y^\a(s)| +|X^a(s)|^2+\Th(s)\,
|X^a(s)|\Big]\, ds, \qq\all t\in [0,T].
\end{equation}
Then using Gronwall's  inequality again, we have
\begin{eqnarray}\label{E507}
\nnb |Y^\a(t)| &\leq & C\int^t_0
\Big[|X^a(s)|^2+\Th(s)\, |X^a(s)|\Big]\, ds\\
\nnb &\leq &  C\int^t_0ds
\int^s_0 \Big\{\Big[\Th(\t)\Big]^2 +\Th(s)\Th(\t) \Big\} \, d\t  \\
\nnb &=& C \int^t_0
(t-s+1)\Big[\Th(s)\Big]^2\, ds\\
&\leq & C \int^t_0 \Big[\Th(s)\Big]^2\,
ds,\qq\all t\in [0,T].
\end{eqnarray}
That is, \eqref{E503} holds.

The proof of \eqref{E504} is similar but a little complex. Consider
\eqref{E404} and \eqref{E406}. We have
\begin{eqnarray}\label{E508}
\nnb &  &   \Big|\int^t_0ds\int^1_0 d\t\int^1_0
 \t \ip{f^k_{xx} (s,\bx(s)+\t\z (x^\a
(s)-\bx(s)),\bu(s))X^\a(s), X^\a(s)} \, d\z \\
\nnb && \qq -{1\over 2}  \int^t_0
   \ip{f^k_{xx} (s,\bx(s) ,\bu(s))X(s), X(s)} \, ds\Big| \\
\nnb & =  &  \Big| \int^t_0ds\int^1_0 d\t\int^1_0 \t \ip{f^k_{xx}
(s,\bx(s)+\t\z (x^\a (s)-\bx(s)),\bu(s))X^\a(s), X^\a(s)} \, d\z\\
\nnb && -\int^t_0ds\int^1_0 d\t\int^1_0 \t \ip{f^k_{xx}
(s,\bx(s),\bu(s))X^\a(s), X^\a(s)} \, d\z\\
 \nnb &    &
+\int^t_0ds\int^1_0 d\t\int^1_0
 \t  \ip{f^k_{xx} (s,\bx(s),\bu(s)) \Big(X^\a(s)+X(s)\Big), \a Y^\a(s)} \, d\z \Big|\\
\nnb & \leq  &   \int^t_0ds\int^1_0 d\t\int^1_0
 \t \,\o\big( \t\z \a |X^\a
(s)|\big)\, |X^\a(s)|^2\, d\z \\
\nnb &    &   +C\int^t_0ds\int^1_0 d\t\int^1_0
 \t \a \Big(|X^\a(s)|+|X(s)|\Big)\, |Y^\a(s)| \, d\z \\
\nnb & \leq  &   \int^t_0 \o\big(|X^\a
(s)|\big)\, |X^\a(s)|^2\, ds  +C\int^t_0 \Big(|X^\a(s)|+|X(s)|\Big)\, |Y^\a(s)| \, ds \\
\nnb & \leq  & C  \int^t_0   \Big[\o\Big(C\int^s_0\Th(\t)\, d\t\Big)+\int^s_0\Th(\t)\, d\t\Big]\, \Big\{\int^s_0\Big[\Th(\t)\Big]^2\, d\t\Big\}\, ds \\
 & \leq  & C    \Big[\o\Big(C\int^T_0\Th(\t)\, d\t\Big)+\int^T_0\Th(\t)\, d\t\Big]\,  \int^t_0\Big[\Th(s)\Big]^2\, ds.
\end{eqnarray}
Similarly,
\begin{eqnarray}\label{E509}
\nnb    & &   \Big|   \int^t_0 ds\int^1_0  \ip{
 f^k_x(s,\bx(s)+\t  (x^\a
(s)-\bx(s)),u(s)), X^\a(s)}\, d\t \\
\nnb  & &  -  \int^t_0 ds\int^1_0 \ip{
 f^k_x(s,\bx(s)+\t  (x^\a
(s)-\bx(s)),\bu(s)),X^\a(s)}  \, d\t\\
\nnb    & &  -  \int^t_0   \ip{
 f^k_x(s,\bx(s),u(s)), X(s)}\, ds   + \int^t_0  \ip{
 f^k_x(s,\bx(s) ,\bu(s)),X (s)}  \, ds\Big| \\
\nnb    & =&   \Big|   \int^t_0 ds\int^1_0 d\t\int^1_0 \ip{
 \a \t f^k_{xx}(s,\bx(s)+\t \z (x^\a
(s)-\bx(s)),u(s))X^\a(s), X^\a(s)}\, d\t \\
\nnb    &  &   -   \int^t_0 ds\int^1_0 d\t\int^1_0 \ip{
 \a\t f^k_{xx}(s,\bx(s)+\t \z (x^\a
(s)-\bx(s)),\bu(s))X^\a(s), X^\a(s)}\, d\t \\
\nnb && + \int^t_0   \ip{
 f^k_x(s,\bx(s),u(s)), \a Y^\a (s)}\, ds  - \int^t_0   \ip{
 f^k_x(s,\bx(s),\bu(s)), \a Y^\a (s)}\, ds\Big|\\
\nnb & \leq  &   \int^t_0 \o\big(|X^\a
(s)|\big)\, |X^\a(s)|^2\, ds  +\int^t_0 \Th(s)\, |Y^\a(s)| \, ds \qq\qq\qq\qq\qq\qq\,\\
\nnb & \leq  & C  \int^t_0   \Big[\o\Big(C\int^s_0\Th(\t)\, d\t\Big)+\Th(s)\Big]\, \Big\{\int^s_0\Big[\Th(\t)\Big]^2\, d\t\Big\}\, ds \\
 & \leq  & C    \Big[\o\Big(C\int^T_0\Th(\t)\, d\t\Big)+\int^T_0\Th(\t)\, d\t\Big]\,  \int^t_0\Big[\Th(s)\Big]^2\, ds.
\end{eqnarray}
While
\begin{eqnarray}\label{E510}
\nnb    & &   \Big| \int^t_0 \ip{f^k_x(s,\bx(s),\bu(s)), Y^\a(s)} \,
ds -\int^t_0 \ip{f^k_x(s,\bx(s),\bu(s)), Y (s)} \, ds\Big|\\
&\leq & \int^t_0 L |Y^\a(s)-Y (s)| \, ds.
\end{eqnarray}
Combining \eqref{E508}---\eqref{E510} with  \eqref{E404} and
\eqref{E406}, we get
\begin{eqnarray}\label{E511}
\nnb   |Y^\a(t)-Y (t)|  & \leq &   C \int^t_0   |Y^\a(s)-Y (s)| \, ds \\
 & & +C \Big[\o\Big(C\int^T_0\Th(\t)\, d\t\Big)+\int^T_0\Th(\t)\, d\t\Big]\,  \int^t_0\Big[\Th(s)\Big]^2\, ds.
\end{eqnarray}
Then Gronwall's inequality implies \eqref{E504}.
\endpf

Now, we give a sufficient optimality condition in the following:

\begin{Theorem}\tt\label{T502}  Assume \thb{S1}---\thb{S4} hold and $\bu(\cd)\in \cU_{ad}$ satisfy \eqref{E305}, \eqref{E307}, \eqref{E308}.  Let $\o(\cd)$ be the function appeared in
 \thb{S4}. If there exists a $\b>0$, such that for any $u(\cd)\in
\cU_{ad}$,
\begin{eqnarray}\label{512}
\nnb  &    & \ds \int^T_0dt \int^t_0\Big\langle
\bW(t)(f(t,\bx(t),\bu(t))-f(t,\bx(t),
 u(t)))\\
 \nnb & & + H_x(t,\bx(t),\bu(t),\bpsi(t))-H_x(t,\bx(t),
u(t),\bpsi(t)),\\
\nnb && \bGP(t)\bGP(s)^{-1}(f(s,\bx(s),\bu(s))-f(s,\bx(s),
 u(s)))\Big\rangle\, ds\\
&\leq &\ds - \b \int^T_0 \Big[\o\big(\rho(u(t),\bu(t))\big)\Big]^2\,
dt,
\end{eqnarray}
then there exists an $\ve_0>0$, such that for any
$$
u(\cd)\in \cV\equiv \set{v(\cd)\in \cU_{ad}\Big|
\int^T_0\o\big(\rho(v(t),\bu(t))\big)\, dt\leq \ve_0 },
$$
\begin{equation}\label{E513}
J(u(\cd))-J(\bu(\cd))\geq {\b\over 2}\int^T_0
\Big[\o\big(\rho(u(t),\bu(t))\big)\Big]^2\, dt.
\end{equation}
In particular, $\bu(\cd)$ minimizes $J(\cd)$ over $\cV$.
\end{Theorem}
\Proof The proof of the theorem is very similar to that of
\eqref{E504} in Lemma \ref{T501}. By \eqref{E308},
\begin{eqnarray}\label{E514}
\nnb & & \int^T_0 \Big[\ip{f^0_x(t,\bx(t),\bu(t)), X(t)}+
 f^0(t,\bx (t),u(t))-f^0(t,\bx (t),\bu(t)) \Big]\, dt\\
 &=& \int^T_0 \Big[H(t,\bx(t),\bu(t)),\bpsi(t))-
H(t,\bx(t),u(t)),\bpsi(t))\Big]\, dt\geq 0.
\end{eqnarray}
Then
\begin{eqnarray}\label{E515}
\nnb  &   & \ds   J((1-\a)\d_{\bu(\cd)}+ \a
 \d_{u(\cd)} )-J( \d_{\bu(\cd)}) \\
 \nnb &=& \int^T_0 \Big[  f^0 (t, x^\a
(t) ,\bu(t))-f^0 (t, \bx
(t) ,\bu(t))+\a \big(f^0(t,x^\a(t),u(t))-f^0(t,x^\a(t),\bu(t))\big) \Big]\, dt\\
\nnb &=&   \a \int^T_0\Big[\int^1_0 \ip{f^0_x(t,\bx(t)+s(x^\a
(t)-\bx(t)),\bu(t)), X^\a(t)}\, ds \\
\nnb & &  \qq +
 f^0(t,x^\a(t),u(t))-f^0(t,x^\a(t),\bu(t)) \Big]\, dt\\
\nnb &\geq &   \a \int^T_0dt\int^1_0  \ip{f^0_x(t,\bx(t)+s(x^\a
(t)-\bx(t)),\bu(t))-f^0_x(t,\bx(t),\bu(t)), X^\a(t)} \, ds \\
\nnb && + \a\int^T_0    \ip{f^0_x(t,\bx(t),\bu(t)), X^\a(t) - X(t)}
\, dt  \\
\nnb && +\a\int^T_0 \big( f^0(t,x^\a(t),u(t))-f^0(t,\bx (t),u(t))\big)\, dt\\
\nnb && -\a\int^T_0 \big( f^0(t,x^\a(t),\bu(t))-f^0(t,\bx (t),\bu(t))\big)\, dt\\
\nnb & = &   \a^2 \int^T_0dt\int^1_0 ds\int^1_0
 s \ip{f^0_{xx} (t,\bx(t)+s\t (x^\a
(t)-\bx(t)),\bu(t))X^\a(t), X^\a(t)} \, d\t \\
\nnb && + \a^2\int^T_0     \ip{f^0_x(t,\bx(t),\bu(t)),
Y^\a(t)} \, dt \\
\nnb & &  +\a^2 \int^T_0 dt\int^1_0  \ip{
 f^0_x(t,\bx(t)+s  (x^\a
(t)-\bx(t)),u(t)), X^\a(t)}\, ds \\
  & &  -\a^2 \int^T_0 dt\int^1_0 \ip{
 f^0_x(t,\bx(t)+s  (x^\a
(t)-\bx(t)),\bu(t)),X^\a(t)}  \, ds,
\end{eqnarray}
By (S2)---(S4) and \eqref{E501}---\eqref{E503},
\begin{eqnarray}\label{E516}
\nnb &  &   \int^T_0dt\int^1_0 ds\int^1_0
 s \ip{f^0_{xx} (t,\bx(t)+s\t (x^\a
(t)-\bx(t)),\bu(t))X^\a(t), X^\a(t)} \, d\t \\
\nnb && -{1\over 2}  \int^T_0
  \ip{f^0_{xx} (t,\bx(t),\bu(t))X (t), X (t)} \, d t \\
\nnb &= & \int^T_0dt\int^1_0 ds\int^1_0
 s \ip{\Big(f^0_{xx} (t,\bx(t)+s\t (x^\a
(t)-\bx(t)),\bu(t))-f^0_{xx} (t,\bx(t),\bu(t))\Big)X^\a(t), X^\a(t)} \, d\t \\
\nnb  & & + \a\int^T_0dt\int^1_0 ds\int^1_0
 s \ip{f^0_{xx} (t,\bx(t),\bu(t))\big(X^\a(t)+X(t)\big), Y^\a(t)} \,
 d\t\\
\nnb &\geq & -C\int^T_0dt\int^1_0 ds\int^1_0
 s \Big[\o \big( s\t \a |X^\a(t)| \big)\, | X^\a(t)|^2 + \big(|X^\a(t)|+|X(t)|\big)\, | Y^\a(t)| \Big]\,
 d\t\\
 \nnb &\geq & -C\int^T_0
 \Big[\o \Big( C \int^t_0  \Th(s)\, ds\Big)+ \int^t_0  \Th(s)\, ds\Big] \, \Big\{\int^t_0 \Big[\Th(s)\Big]^2\,
ds\Big\}\,
 dt\\
&\geq & -C
 \Big[\o \Big( C \int^T_0  \Th(s)\, ds\Big)+ \int^T_0  \Th(s)\, ds\Big] \, \Big\{\int^T_0 \Big[\Th(s)\Big]^2\,
ds.
\end{eqnarray}
By (S2)---(S4) and  \eqref{E504},
\begin{eqnarray}\label{E517}
\nnb &   &   \int^T_0     \ip{f^0_x(t,\bx(t),\bu(t)), Y^\a(t)} \, dt
-\int^T_0     \ip{f^0_x(t,\bx(t),\bu(t)),
Y (t)} \, dt\\
\nnb &\geq  & -C  \int^T_0 |Y^\a(t)-Y(t)|\, dt\\
&\geq & - C\Big[\o\Big(C\int^T_0 \Th(s)\, ds \Big)   +   \int^T_0 \Th(s)\, ds \Big]\,
\int^T_0 \Big[\Th(s)\Big]^2\, ds.
\end{eqnarray}
On the other hand, it follows from (S2)---(S4) and \eqref{E503}---\eqref{E504} that
\begin{eqnarray}\label{E518}
\nnb &   &   \int^T_0 dt\int^1_0  \ip{
 f^0_x(t,\bx(t)+s  (x^\a
(t)-\bx(t)),u(t)), X^\a(t)}\, ds \\
\nnb  & &  -  \int^T_0 dt\int^1_0 \ip{
 f^0_x(t,\bx(t)+s  (x^\a
(t)-\bx(t)),\bu(t)),X^\a(t)}  \, ds\\
\nnb & & - \int^T_0    \ip{
 f^0_x(t,\bx(t),u(t)), X (t)}\, dt   + \int^T_0   \ip{
 f^0_x(t,\bx(t),\bu(t)),X (t)}  \, dt\\
\nnb &=& \int^T_0 dt\int^1_0 ds\int^1_0 s\a \ip{
 f^0_{xx}(t,\bx(t)+s \t (x^\a
(t)-\bx(t)),u(t))X^\a(t), X^\a(t)}\, d\t \\
\nnb & & -\int^T_0 dt\int^1_0 ds\int^1_0 s\a \ip{
 f^0_{xx}(t,\bx(t)+s \t (x^\a
(t)-\bx(t)),\bu(t))X^\a(t), X^\a(t)}\, d\t \\
\nnb  & &  +  \int^T_0 dt\int^1_0 \ip{f^0_x(t,\bx(t),u(t))-
 f^0_x(t,\bx(t),\bu(t)),\a Y^\a(t)}  \, ds\\
\nnb &\geq & -C\int^T_0 \Th(t)\,\Big[|X^\a(t)|^2+|Y^\a(t)\Big]\, dt \\
\nnb &\geq  &-C\int^T_0 \Th(t)\,\Big\{\int^t_0 \Big[\Th(s)\Big]^2\, ds  \Big\}\, dt \\
  & \geq & -C\int^T_0 \Th(t)\,dt \, \int^T_0 \Big[\Th(s)\Big]^2\, ds .
 \end{eqnarray}
Therefore, combining \eqref{E515}---\eqref{E518} with \eqref{E511},
\eqref{E407} and \eqref{E410}, we have
\begin{eqnarray}\label{E519}
\nnb  &   & \ds   {J((1-\a)\d_{\bu(\cd)}+ \a
 \d_{u(\cd)} )-J( \d_{\bu(\cd)})\over \a^2}\\
\nnb &\geq & \int^T_0     \ip{f^0_x(t,\bx(t),\bu(t)), Y (t)} \, dt    +
\int^T_0   \ip{
 f^0_x(t,\bx(t),u(t))- f^0_x(t,\bx(t) ,\bu(t)), X (t)}\, dt\\
\nnb & & +{1\over 2}  \int^T_0
  \ip{f^0_{xx} (t,\bx(t),\bu(t))X (t), X (t)} \, d t\\
\nnb && -C\Big[\o\Big(C\int^T_0 \Th(t)\, dt \Big)   +   \int^T_0 \Th(t)\, dt \Big]\,
\int^T_0 \Big[\Th(t)\Big]^2\, ds\\
\nnb &=& \ds \int^T_0dt \int^t_0\Big\langle
\bW(t)(f(t,\bx(t),\bu(t))-f(t,\bx(t),
 u(t)))+ H_x(t,\bx(t),\bu(t),\bpsi(t))\\
 \nnb & & -H_x(t,\bx(t),
u(t),\bpsi(t)), \bGP(t)\bGP(s)^{-1}
 (f(s,\bx(s),u(s))-f(s,\bx(s),
 \bu(s)))\Big\rangle\, ds\\
\nnb && -C\Big[\o\Big(C\int^T_0 \Th(t)\, dt \Big)   +   \int^T_0 \Th(t)\, dt \Big]\,\int^T_0 \Big[\Th(t)\Big]^2\, dt\\
  &\geq & \Big[\b-C\,\o\Big(C\int^T_0 \Th(t)\, dt \Big) -C \int^T_0 \Th(t)\, dt \Big]\, \int^T_0 \Big[\Th(t)\Big]^2\, dt.
 \end{eqnarray}
Since
$$
\lim_{\ve\to 0^+} C\Big[ \o (C\ve  ) +\ve \Big]=0,
$$
there exists an $\ve_0>0$, independent of $\a\in (0,1]$ and $u(\cd)$,  such that when
$$
\int^T_0\Th(t)\, dt\leq \ve_0,
$$
it holds that
\begin{eqnarray}\label{E520}
  &   & \ds   {J((1-\a)\d_{\bu(\cd)}+ \a
 \d_{u(\cd)} )-J( \d_{\bu(\cd)})\over \a^2} \geq {\b\over 2}\, \int^T_0 \Big[\Th(t)\Big]^2\, dt.
\end{eqnarray}
Choosing $\a=1$ in \eqref{E520}, we get \eqref{E513} and finish the
proof.
\endpf

 \footnotesize
\vspace{5mm}\footnotesize \ \\


\begin{thebibliography}{99}
\bibitem{Bon-Her1} J. F.  Bonnans and A. Hermant, No-gap second-order optimality conditions for optimal
control problems with a single state constraint and control, Math. Program., Ser. B, 117(2009), pp. 21--50.
\bibitem{Bon-Her2} J. F.  Bonnans and A. Hermant, Second-order analysis for optimal control problems
with pure state constraints and mixed control-state constraints, Ann. I. H. Poincare,  26(2009), pp. 561--598.
\bibitem{Cas-Car-Tr} E. Casas, J.C. de Los Reye and F. Tr\"oltzsch, Sufficient second-order optimality conditions for
semilinear control problems with pointwise state
constraints, SIAM J. Optim., 19 (2008), pp. 616--643.
\bibitem{Cas-Mat} E. Casas and M. Mateos, Second order optimality conditions for semilinear elliptic control
problems with finitely many state constraints, SIAM J. Control Optim., 40 (2002), pp. 1431--1454.
\bibitem{Cas-Tr1}
 E. Casas and F. Tr\"oltzsch, Second order necessary optimality conditions for some state-
constrained control problems of semilinear elliptic equations, Appl. Math. Optim., 39 (1999), pp.
211--227.
\bibitem{Cas-Tr2} E. Casas and F. Tr\"oltzsch, Second-order necessary and sufficient optimality conditions for opti-
mality conditions for optimization problems and applications to control theory, SIAM J. Optim.,
13 (2002), pp. 406--431.
\bibitem{Cas-Tr3} E. Casas and F. Tr\"oltzsch, First- and second-order optimality conditions for a
class of optimal control problems with quasilinear elliptic equations, SIAM J. Control Optim., 48 (2009), pp. 688--718.
\bibitem{Cas-Tr-Un} E. Casas, F. Tr\"oltzsch  and A. Unger, Second-order sufficientt optimality conditions for some
state-constrained control problems of semilinear elliptic equations, SIAM J. Control Optim., 38
(2000), pp. 1369--1391.
\bibitem{Fa1}  H. O. Fattorini,  Relaxed controls in infinite
     dimensional systems, International Series of Numerical
     Mathematics, 100 (1991), pp.  115--128.
\bibitem{Ga} R. Gamkrelidze,  ``Principle of Optimal Control
          Theory", Plenum Press,     New York, 1978.
\bibitem{G-K} R. Gabasov and F. M. Kirillova, High order necessary conditions for pptimality, SIAM J. Control, 10 (1972), pp. 127--168.
\bibitem{Hu} D. G. Hull, Sufficient conditions for a minimum of the free-final-time optimal control problem,
J. Optim. Theory Appl., 68 (1991), pp. 275--287.
\bibitem{Ke}
H. J. Kelly, A second variation test for singular extremals, AIAA J., 2 (1964), pp. 1380--1382.
\bibitem{Ko-Mo} R. E. Kopp and H. G. Moyer, Necessary conditions for singular extremals, AIAA J., 3 (1965), pp.
1439--1444.
\bibitem{Kre} A. J. Krener, The high order maximal principle and its application to singular extremals, SIAM
J. Control Optim., 15 (1977), pp. 256--293.
\bibitem{Li-Yong} X. Li  and J. Yong, ``Optimal Control Theory for Infinite
     Dimensional Systems", Birkh\"{a}user, Boston, 1995.
\bibitem{Lou-Yong 2009} H. Lou  and J. Yong,  Optimality Conditions for
Semilinear  Elliptic Equations with Leading Term Containing
Controls,   SIAM J. Control Optim.,   48 (2009),  pp. 2366--2387 .
\bibitem{Mal} K. Malanowski, Sufficient optimality conditions for optimal control subject to state constraints,
SIAM J. Control Optim., 35 (1997), pp. 205--227.
\bibitem{Mal-Mau-Pi} K. Malanowski, H. Maurer and S. Pickenhain, Second-order sufficient conditions for state-constrained optimal control problems,
J. Optim. Theory Appl., 123(2004), pp. 595--617.
\bibitem{Mau-Ob} H. Maurer and H. J. Oberle, Second order su¡Àcient conditions for optimal control problems
with free final time: the Riccati approach, SIAM J. Control Optim., 41(2002), pp. 380--403.
\bibitem{Mau-Pes} H. Maurer and H. J. Pesch, Solution differentiability for parametric nonlinear control problems
with control-state constraints, J. Optim. Theory Appl., 86 (1995), pp. 285--309.
\bibitem{Mau-Pick} H. Maurer and S. Pickenhain, Second order sufficient conditions for optimal control problems
with mixed control-state constraints, J. Optim. Theory Appl., 86 (1995), pp. 649--667.
\bibitem{Mau-Os} H. Maurer and N. P. Osmolovskii, Second order sufficient conditions for time-optimal bang-
bang control, SIAM J. Optim., 42 (2004), pp. 2239--2263.
\bibitem{Pal-Zei}
Z. Pales and V. Zeidan,
 The critical tangent cone in second-order conditions for optimal control,
Nonlinear Analysis, 47(2001), pp. 1149-1161.

\bibitem{Pick} S. Pickenhain, Sufficiency conditions for weak local minima in multidimensional optimal control
problems with mixed control-state restrictions, Z. Anal. Anwendungen, 11 (1992), pp.
559--568.
\bibitem{Ray-Tr} J. P. Raymond  and F. Trltzsch, Second order su¡Àcient optimality conditions for nonlinear
parabolic control problems with state constraints, Discrete Contin. Dynam. Systems, 6 (2000), pp.
431--450.
\bibitem{Ro-Tr1}
A. R\"osch and F. Tr\"oltzsch, Sufficient second-order optimality conditions for a
parabolic optimal control problem with pointwise control-state constraints, SIAM J. Control Optim., 42 (2003), pp. 138--154.
\bibitem{Ro-Tr2} A. R\"osch and F. Tr\"oltzsch,  Sufficient second-order optimality conditions for an elliptic optimal control problem with pointwise
control-state constraints, SIAM J.  Optim., 17 (2006), pp. 776--794.
\bibitem{Ros}
J. F. Rosenblueth, A new derivation of second-order conditions for equality
control constraints, Appl. Math. Letters,  21 (2008), pp. 910--915.
\bibitem{Wach} D.  Wachsmuth,
Sufficient second-order optimality conditions for convex control constraints, J. Math. Anal. Appl.,  319(2006), pp. 228--247.
\bibitem{Wa3}  J. Warga,   ``Optimal Control  of Differential and
Functional Equations",   Academic Press, New York,  1972.
\bibitem{Wang-He} L. Wang and P. He, Second-order optimality conditions for
optimal control problems governed by
3-dimensional Nevier-Stokes equations, Acta Math. Scientia,   26B(2006), pp. 729--734.

\bibitem{Zygmund} A. Zygmund, ``Trigonometric Series", 3rd ed.,  Cambridge University Press, Cambridge, 2002.




\end{thebibliography}
\end{document}